\documentclass[12pt]{article}[1994/06/01]
\usepackage{amsxtra,amssymb,amsthm,amsmath,latexsym}
\usepackage{graphicx}
\newtheorem{thm}{Theorem}[section]

\newcommand{\R}{{\mathbb R}}

\newcommand{\bee}{\begin{equation*}}
\newcommand{\eee}{\end{equation*}}
\newcommand{\be}{\begin{equation}}
\newcommand{\ee}{\end{equation}}
\begin{document}
\title{A DSM proof of surjectivity of monotone nonlinear mappings}
\author{A. G. Ramm}
\date{}
\maketitle
\section{Introduction}
It is well-known that a continuous monotone function $f:\R\to \R$,
such that \be\label{e1} \lim_{|x|\to \infty} xf(x)/|x|=\infty, \ee is
surjective, i.e., the equation $f(x)=y$ is solvable for any $y\in
\R$. Indeed, the monotonicity of $f$ implies \be\label{e2}
[f(x)-f(s)](x-s)\geq 0,\quad \forall x,s\in \R. \ee Therefore,
taking $y=0$ without loss of generality, one concludes from
(\ref{e1})  that $f(x)\leq 0$ for $x\leq 0$ and $f(x)\geq
0$ for $x\geq 0$. Since $f$ is continuous, it follows that there is
a point $x_0$ such that $f(x_0)=0.$\\
If $y\neq 0$ is an arbitrary real number, then the function
$F(x)=f(x)-y$ satisfies inequality (\ref{e2}) with $F$ in place of
$f$, provided that (\ref{e2}) holds for $f$. Condition (\ref{e1}) is 
also satisfied 
for $F$ if it holds for $f$:
$$\lim_{|x|\to
\infty}\frac{xF(x)}{|x|}=\lim_{|x|\to
\infty}\left(\frac{xf(x)}{|x|}-\frac{xy}{|x|}\right)=\infty.$$
Conditions (\ref{e1}) and (\ref{e2}) are generalized for nonlinear
mappings $F$ in a real Hilbert space $H$ as follows: \be\label{e3}
\lim_{\|u\|\to \infty}\frac{(u,F(u))}{\|u\|}=\infty, \ee and
\be\label{e4} (F(u)-F(v),u-v)\geq 0\quad \forall u,v\in H. \ee 
here $(u,v)$ stands for the inner product in $H$.

We
want to prove that if $F$ is twice Frechet differentiable and
conditions (\ref{e3})-(\ref{e4}) hold, then $F$ is surjective, i.e.,
the equation \be\label{e5} F(u)=h\ee is solvable for every $h\in H$.
This is a basic result in the theory of monotone operators (see,
e.g., \cite{D}, \cite{Z}). Our aim is to give a simple and short
proof of this result based on the Dynamical Systems Method (DSM)
developed in \cite{R499}.
\begin{thm}
Assume that $F:H\to H$ is Fr$\acute{e}$chet differentiable mapping
satisfying conditions (\ref{e3}), (\ref{e4}). Then equation
(\ref{e5}) is solvable for any $h$.
\end{thm}
{bf Remark 1.} If in (\ref{e4}) one has a strict inequality for $u\neq
v$, then the solution to (\ref{e5}) is unique.

{\bf Remark 2.} Condition (\ref{e4}) and Fr$\acute{e}$chet
differentiability imply that $A:=F'(u)\geq 0\quad \forall u\in H$

{\bf Remark 3.} The Fr$\acute{e}$chet differentiability assumption can be
weakened to semi-continuity (see, e.g., \cite{D}), but then the
proof loses its elementary character.

Our proof is elementary and is suitable for
undergraduate students: it requires a very limited background.

\section{Proof}
Let us formulate the steps of our proof.\\
{\t Step 1.} For any $a=const>0$ the equation
\be\label{e6}F(u_a)+au_a=h\ee has a unique solution $u_a.$\\
{\it Step 2.} \be\label{e7} \sup_{0<a<1}\|u_a\|<c,\quad c=const>0.\ee By
$c$ we denote various constants independent of $a$.\\
{\it Step 3.} Using (\ref{e7}), select a sequence $u_n=u_{a_n},\ a_n\to
0,$ weakly convergent in $H$  to an element
$u$:
\be\label{e8}u_n\rightharpoonup u,\qquad n\to \infty.\ee 
From
(\ref{e6}) and (\ref{e8}) it follows that 
\be\label{e9} F(u_n)\to
h,\qquad n\to \infty.\ee From (\ref{e8}), (\ref{e9}) and (\ref{e4})
one concludes that $u$ solves (\ref{e5}).\\
Let us give a {\it detailed proof}.\\
{\it Step 1.} Consider the problem: 
\be\label{e10}
\dot{v}=-A_a^{-1}[F(v)+av-h],\quad v(0)=0.\ee 
Here $\dot{v}:=\frac{dv}{dt}$, $A_a:=A+aI,\
A:=F'(v).$ Problem (\ref{e10}) is a version of a DSM. We claim that:
\begin{itemize}
\item[a)] problem (\ref{e10}) has a unique global solution , i.e. ,
the solution defined for all $t\in [0,\infty)$,
\item[b)] there exists $v(\infty):=\lim_{t\to \infty}v(t)$, and
\item[c)] $F(v(\infty))+av(\infty)=h.$
\end{itemize} Claim $a)$ follows from local solvability of problem
(\ref{e10}) and a uniform with respect to $t$ bound on the norm
$\|v(t)\|.$ This bound is obtained below (see (\ref{e13})).

Denote $\|F(v(t))+av(t)-h\|:=g(t),\ \dot{g}:=\frac{dg}{dt}.$ Using
(\ref{e10}), one gets
$$g\dot{g}=\left((F'(v)+aI)\dot{v},F(v(t))+av(t)-h\right)=-g^2.$$
Thus \be\label{e11} g(t)=g(0)e^{-t}. \ee From (\ref{e11}) and
(\ref{e10}) one gets \be\label{e12} \|\dot{v}\|\leq
\frac{g(0)}{a}e^{-t},\ee where the estimate $\|A_a^{-1}\|\leq
\frac{1}{a}$ was used. This estimate holds because $A=F'(v(t))\geq
0$ by the monotonicity of $F$. Integrating (\ref{e12}) from $t$ to
infinity yields \be\label{e13} \|v(t)-v(\infty)\|\leq
\frac{g(0)}{a}e^{-t}. \ee Note that if $\|\dot{v}\|\leq g(t)$ and
$g(t)\in L^1(0,\infty)$, then $v(\infty)$ exists by the Cauchy
criterion for the existence of a limit:
$$\|v(t)-v(s)\|\leq
\int_s^tg(\tau) d\tau\to 0,\ t,s\to \infty, \ t>s.$$ 
It follows from
(\ref{e12}) that \be\label{e14} \lim_{t\to \infty}\|\dot{v}\|=0. \ee
Therefore, passing to the limit $t\to \infty$ in (\ref{e10}), one
gets \be\label{e15}
0=-A_a^{-1}(v(\infty))[F(v(\infty))+av(\infty)-h].\ee Applying the
operator $A_a(v(\infty))$ to equation (\ref{e15}), one sees that
$v(\infty)$ solves equation (\ref{e6}). Uniqueness of the solution
to (\ref{e6}) is easy to prove: if $v$ and $w$ solve (\ref{e6}),
then $$F(v)-F(w)+a(v-w)=0,\ a>0.$$ 
Multiply this equation by $v-w$,
use the monotonicity of $F$ (see (\ref{e4})), and conclude that $v=w$. Step 1 
is completed.\\
{\it Step 2.} Let us prove (\ref{e7}). Multiply (\ref{e6}) by
$\frac{u_a}{\|u_a\|}$ and get \be\label{e16}
\frac{(F(u_a),u_a)}{\|u_a\|}+a\|u_a\|=\frac{(f,u_a)}{\|u_a\|}. \ee
Since $a>0$ and $\frac{(f,u_a)}{\|u_a\|}\leq \|f\|$, one gets
\be\label{e17} \frac{(F(u_a),u_a)}{\|u_a\|}\leq \|f\|. \ee From
(\ref{e17}) and (\ref{e3}) the desired estimate (\ref{e7}) follows.
Step 2 is completed. \\
{\it Step 3.} Let us prove that (\ref{e4}), (\ref{e8}) and (\ref{e9})
imply (\ref{e5}). Let $\eta\in H$ be arbitrary, and $s>0$ be a small
number. Note that $u_n\rightharpoonup u$ and $g_n\to g$ imply
$(u_n,g_n)\to (u,g).$ Using (\ref{e4}), one gets: \be\label{e18}
\left(F(u_n)-F(u-s\eta),u_n-u+s\eta\right)\geq 0,\quad \forall
\eta\in H,\ s>0 .\ee Let $n\to \infty$ in (\ref{e18}). Then, using
(\ref{e8}) and (\ref{e9}), one concludes that
\be\label{e19}\begin{split} &(h-F(u-s\eta),s\eta)\geq 0\quad \forall
\eta\in H,\ s>0,\\
&\text{or}\\
&(h-F(u-s\eta),\eta)\geq 0\quad \forall
\eta\in H,\ s>0.\\
\end{split}
\ee Let $s\to 0$ and use the continuity of $F$. Then (\ref{e19})
implies \be\label{e20} (h-F(u),\eta)\geq 0\quad \forall \eta\in
H.\ee Taking $\eta=h-F(u)$ in (\ref{e20}), one concludes that
$F(u)=h$. Step 3 is completed. Theorem 1 is proved. $\Box$

\end{document}